\def\ARXIV{ARXIV}
\def\JOURNAL{JOURNAL}
\def\form{ARXIV}

\ifx\form\ARXIV
\documentclass[12pt]{amsart}
\usepackage{graphicx,setspace,algorithmic,color,amssymb,hyperref}
\fi
\ifx\form\JOURNAL
\documentclass{ws-procs9x6}
\usepackage{graphicx,setspace,algorithmic,color,amssymb,hyperref}
\fi

\newenvironment{M2}{ \begin{spacing}{0.4} 
\bigskip \small} { \smallskip 
\end{spacing} }

\ifx\form\ARXIV
\newtheorem{theorem}{Theorem}
\theoremstyle{plain}

\newtheorem{algorithm}{Algorithm}

\newtheorem{definition}{Definition}
\newtheorem{example}{Example}

\newtheorem{experiment}{Experiment}

\newtheorem{remark}{Remark}

\numberwithin{equation}{section}
\fi
\ifx\form\JOURNAL
\newtheorem{experiment}{Experiment}
\renewcommand{\cite}[2][]{{[\ifthenelse{\equal{#1}{}}{}{#1,~}\refcite{#2}]}}
\fi

\newcommand{\bem}[1]{{\em {#1}}}

\def \gets{\leftarrow}
\def \<{\langle}
\def \>{\rangle}

\def\Ann{\operatorname{Ann}}

\def\WeylCl{\operatorname{WeylCl}}

\def\gr{\operatorname{gr}}

\def\ord{\operatorname{ord}}

\newcommand\KappaAndAnnihilator{{\sc KappaAndAnnihilator}}
\newcommand\CheckGenericity{{\sc CheckGenericity}}

\newcommand{\bC}{{\mathbb C}}
\newcommand{\bN}{{\mathbb N}}

\newcommand{\ZZ}{{\mathbb Z}}
\newcommand{\RR}{{\mathbb R}}
\newcommand{\bQ}{{\mathbb Q}}

\newcommand{\bV}{{\mathbb V}}

\newcommand{\cO}{{\mathcal  O}}
\newcommand{\cD}{{\mathcal D}}

\newcommand{\bfp}{{\boldsymbol{\partial}}}

\newcommand{\bfx}{{\boldsymbol{x}}}

\newcommand{\bfxi}{{\boldsymbol{\xi}}}

\newcommand{\zero}{{\mathbf{0}}}

\newcommand{\p}{{\partial}}

\newcommand{\IN}{{\operatorname{in}}}

 \def\D0{D_\zero}
\def\Span{\operatorname{Span}}

\newcommand{\fa}{{\mathfrak a}}

\newcommand{\calI}{{\mathcal I}}
\newcommand{\calA}{{\mathcal A}}

\newcommand{\Dx}{{D_X}}

\def\cI{\gr_w}
\def\cV{\operatorname{Char}}
\def\CChar{\operatorname{CChar}}

\def\D4M2{{\tt Dmodules}}

\catcode`\@=11\relax
\newdimen\p@renwd
\font\tenex=cmex10 \setbox0=\hbox{\tenex B} \p@renwd=\wd0
\def\bbordermatrix#1{\begingroup \m@th
\setbox\z@\vbox{\def\\{\crcr\noalign{\kern2\p@\global\let\cr\endline}}%
    \ialign{$##$\hfil\kern2\p@\kern\p@renwd&\thinspace\hfil$##$\hfil
      &&\quad\hfil$##$\hfil\crcr
      \omit\strut\hfil\crcr\noalign{\kern-\baselineskip}%
      #1\crcr\omit\strut\cr}}%
  \setbox\tw@\vbox{\unvcopy\z@\global\setbox\@ne\lastbox}%
  \setbox\tw@\hbox{\unhbox\@ne\unskip\global\setbox\@ne\lastbox}%
  \setbox\tw@\hbox{$\kern\wd\@ne\kern-\p@renwd\left[\kern-\wd\@ne
    \global\setbox\@ne\vbox{\box\@ne\kern2\p@}%
    \vcenter{\kern-\ht\@ne\unvbox\z@\kern-\baselineskip}\,\right]$}%
  \null\;\vbox{\kern\ht\@ne\box\tw@}\endgroup}
\catcode`\@=12\relax

\newlength\cellsize 
\savebox2{%
\begin{picture}(10,10)
\put(0,0){\line(1,0){10}}
\put(0,0){\line(0,1){10}}
\put(10,0){\line(0,1){10}}
\put(0,10){\line(1,0){10}}
\end{picture}}
\newcommand\cellify[1]{\def\thearg{#1}\def\nothing{}%
\ifx\thearg\nothing
\vrule width0pt height\cellsize depth0pt\else
\hbox to 0pt{\usebox2\hss}\fi%
\vbox to 10\unitlength{
\vss
\hbox to 10\unitlength{\hss$#1$\hss}
\vss}}
\newcommand\tableau[1]{\vtop{\let\\=\cr
\setlength\baselineskip{-16000pt}
\setlength\lineskiplimit{16000pt}
\setlength\lineskip{0pt}
\halign{&\cellify{##}\cr#1\crcr}}}
\savebox3{%
\begin{picture}(15,15)
\put(0,0){\line(1,0){15}}
\put(0,0){\line(0,1){15}}
\put(15,0){\line(0,1){15}}
\put(0,15){\line(1,0){15}}
\end{picture}}
\newcommand\expath[1]{%
\hbox to 0pt{\usebox3\hss}%
\vbox to 15\unitlength{
\vss
\hbox to 15\unitlength{\hss$#1$\hss}
\vss}}

\begin{document}

\title {Computing localizations iteratively}

\author{Francisco-Jes\'us Castro-Jim\'enez}
\ifx\form\ARXIV
\thanks{F.-J. Castro-Jim\'{e}nez. Depto. \'{A}lgebra. Universidad de Sevilla, Sevilla, Spain({\tt castro@us.es}).
         Partially supported by MTM2010-19336 and FEDER, FQM-333 and FQM-5849 Junta de Andaluc\'{\i}a}
\fi
\ifx\form\JOURNAL
\address{
Depto. \'{A}lgebra\\ Universidad de Sevilla \\ Sevilla, Spain \\ \tt castro@us.es
}
\fi

\author{Anton Leykin}
\ifx\form\ARXIV
        \thanks{Anton Leykin. School of Mathematics, Georgia Tech, Atlanta GA, USA ({\tt leykin@math.gatech.edu}).
         Partially supported by NSF grant DMS-0914802}
\fi
\ifx\form\JOURNAL
\address{School of Mathematics\\ Georgia Tech\\ Atlanta GA, USA\\ \tt leykin@math.gatech.edu}
\fi

\ifx\form\ARXIV
\maketitle
\fi

\ifx\form\JOURNAL
\bodymatter
\fi

\begin{abstract}
Let $R=\bC[\bfx]$ be a polynomial ring with complex coefficients and $\Dx = \bC\<\bfx,\bfp\>$ be the Weyl algebra. Describing the localization $R_f = R[f^{-1}]$ for nonzero $f\in R$ as a $\Dx$-module amounts to computing the annihilator $A = \Ann(f^a)\subset \Dx$ of the cyclic generator $f^{a}$ for a suitable negative integer $a$. We construct an iterative algorithm that uses truncated annihilators to build $A$ for planar curves.
\end{abstract}

\section*{Introduction}
Let $f\in R:=\bC[x_1,\ldots,x_n]$ be a nonzero polynomial ($n$ being a positive integer) and $R_f$ the localization of the polynomial ring $R$ with respect to $f$. Elements in $R_f$ are rational functions $\frac{g}{f^k}$ with $g\in R$ and $k\in \bN$. If $f$ is not a constant, the $R$--module $R_f$ is not finitely generated. One fundamental result by J. Bernstein states that $R_f$ is finitely generated when considered as a left module over the complex Weyl algebra $A_n:=A_n(\bC)$ of order $n$ (see the needed definitions and precise statements in the next Section). This is a
fundamental result in $\mathcal D$--module theory (the algebraic theory of systems of linear partial differential equations) which can be considered as part of Singularity Theory.

Bernstein's result states even more: the module $R_f$ is cyclic over $A_n$ and in fact it is generated by a rational function of the form $f^{a}$ for some negative integer $a$. That means that $R_f$ can be described as a quotient $\frac{A_n}{Ann(f^a)}$ where $Ann(f^a):=Ann_{A_n}(f^a)$ is the left ideal in $A_n$ formed of the linear differential operators annihilating the rational function $f^a$. By \cite{Kashiwara-78} (see also \cite{Varchenko-82} and \cite{Saito-94}) one can take $a= -n+1$.

One main problem in algorithmic $\mathcal D$--module theory  is to compute a finite system of generators of the annihilating ideal $Ann(f^a)$. There are several algorithms solving this problem (see e.g. \cite{Oaku:local-b}, \cite{Oaku-Takayama-Walther-2000}, \cite{Oaku-Takayama-2001}), which use Gr\"obner bases and elimination theory in the Weyl algebra $A_n$.  In the worst case, computing Gr\"obner bases has a doubly exponential complexity in both (commutative) polynomial rings and the Weyl algebra (see \cite{AschLeykin:GBB} and \cite{Grigoriev-Chistov:D-complexity}); however, in practice, one should expect much longer running times for the latter on the input of the same size.

As the Weyl algebra $A_n$ is a Noetherian ring, one can associate to the couple  $(f,a)$ the smallest integer $\kappa(f^a)$ such that $Ann(f^a)$ is generated by operators of order less than or equal to $\kappa(f^a)$. This numerical invariant plays a relevant role in the so called Logarithmic Comparison Problem with respect to the hypersurface defined in the complex affine space $\bC^n$ by the polynomial equation $f=0$ (see e.g. \cite{Castro-Ucha-Experimental-2004}, \cite{Torrelli-2004}).

In this paper we describe a new algorithm computing $Ann(f^{-1})$ for any reduced complex polynomial $f=f(x,y)$ in two variables. This algorithm uses Gr\"obner bases techniques in both polynomial rings  and the Weyl algebra $A_2$, but avoids elimination theory in the latter: our experiments show that the bottleneck of the algorithm is a syzygy module computation over the former.

We first compute iteratively, for $d=1,2,\ldots$\,, truncated annihilating ideals $Ann^{(d)} (f^{-1})$ generated by linear differential operators annihilating $f^{-1}$ and of order less than or equal to $d$. The ideal $Ann^{(d)} (f^{-1})$ can be computed by using polynomial Gr\"obner basis. We then compute equations for the characteristic cycle of the ideal $\Ann^{(d)} (f^{-1})$ and compare it with the one of $Ann^{} (f^{-1})$ which only depends on the multiplicity of the plane curve $f=0$ at its singular points; for the beforementioned  comparison we only need to localize at each of these singular points.

The paper is organized as follows. In the first Section we survey basic results on the Weyl algebra, holonomic $D$-modules, $b$-functions (or Bernstein-Sato polynomials) and the annihilating ideal of some rational functions on $\bC^n$. In Section 2 we describe the new iterative algorithm, we give a stopping criterion and prove its correctness. We also treat some examples, in particular the so called family of Reiffen's curves $f_{p,q}$. Our implementation of the algorithm in Macaulay2 works on larger examples of Reiffen's curve family than the known general algorithm and our experiment demonstrates an interesting (very simple!) dependence of the order of generation of $\Ann^{(d)}(f_{p,q}^{-1})$ on $p$. In the last Section we conclude the discussion and propose some open questions on the subject.

\section{Preliminaries}

\subsection{Weyl algebra}

Let $n\geq 0$ be an integer and   $X = \bC^n$ be the $n$-dimensional complex affine space.
Define the $n$-th \bem{Weyl algebra} as the associative algebra
$$\Dx = \bC\<\bfx,\bfp\> = \bC\< x_1,\ldots,x_n,\p_1,\ldots,\p_n\>$$
where  $[\partial_i, x_i] = \partial_ix_i-x_i\partial_i=1$ and all other pairs of generators commute.
The Weyl algebra $\Dx$ is isomorphic to the algebra of \bem{linear differential operators} with coefficients in the polynomial ring $R=\bC[\bfx]=\bC
[x_1,\ldots,x_n]$.
Every element in $\Dx$ has a unique \bem{normal form}
  $$
  Q = \sum_{\alpha,\beta\in\bN^n} q_{\alpha\beta} \bfx^\alpha\bfp^\beta,
  $$
where finitely many of $q_{\alpha\beta}\in \bC$ are nonzero and where $\bfx^\alpha\bfp^\beta$ stands for $x_1^{\alpha_1}\cdots x_n^{\alpha_n} \p_1^{\beta_1} \cdots \p_n^{\beta_n}$.
We denote the $n$-th Weyl algebra $A_n$ if we need to emphasize the dimension of $X$; throughout this paper $\Dx = A_2$ is used the most.

The ring $\Dx$ is simple: there are only trivial two-sided ideals. 
All $\Dx$-ideals and $\Dx$-modules considered in this article are \bem{left} ideals and modules, respectively.

Examples of $\Dx$-modules include functional spaces, e.g., polynomial functions $\bC[\bfx]$ and smooth functions $C^\infty(X)$, and power series rings $\bC[[\bfx]]$.

Another example is the localization of the polynomial ring $\bC[\bfx,f^{-1}]$ where $f$ is a nonzero polynomial with the natural action defined as follows:
\begin{align*} x_i\cdot g f^{-j} &= x_i g f^{-j},\\
\p_i \cdot g f^{-j} &= \left(\frac{\p g}{\p x_i}f - j g \frac{\p f}{\p x_i}\right) {f^{-j-1}},
\end{align*}
for $1\leq i \leq n$, $g\in \bC[\bfx]$, and $j \in \bN$.

\subsection{Gr\"obner bases}

It is possible to compute in the Weyl algebra $\Dx$, since it is \bem{Gr\"obner friendly}; see e.g. \cite{Castro_thesis, Castro87}, \cite{SST} and \cite{Levandovskyy-PhD-2005}.

Gr\"obner bases can be computed with respect to any $w$-compatible monomial order, where $w=(w_\bfx,w_\bfp)\in \RR^{2n}$ satisfies $w_\bfx+w_\bfp \geq \zero$ componentwise. This condition, in particular, guarantees that the filtration $\{F^w_{i}\}_{i\in \ZZ}$ of $\Dx$, where
$$
F^w_{i} = \Span\left\{ \bfx^\alpha\bfp^\beta : \alpha,\beta\in\bN^n, \ w\cdot(\alpha,\beta)\leq i \right\},
$$
is preserved under taking the commutator. Therefore, the associated graded algebra $\gr_w(\Dx)$ is well defined. Note that
\begin{itemize}
\item if $w_\bfx+w_\bfp = \zero$ then $\gr_w(\Dx) = \Dx$;
\item if $w_\bfx+w_\bfp > \zero$ componentwise then $$\gr_w(\Dx) \cong \bC[\bfx,\bfxi] = \bC[x_1,\ldots,x_n,\xi_1,\ldots,\xi_n],$$ a polynomial ring in $2n$ variables.
\end{itemize}

The $w$-order of  a nonzero operator $Q= \sum_{\alpha,\beta\in\bN^n} q_{\alpha\beta} \bfx^\alpha\bfp^\beta \in \Dx$, denoted $\ord_w(Q)$,  is defined as the maximum of $w\cdot (\alpha,\beta)$ for $q_{\alpha\beta}\not=0$. The initial part of $Q$ with respect to $w$ is $$\IN_w(Q)= \sum_{w\cdot (\alpha,\beta) = \ord_w(Q)} q_{\alpha\beta} \bfx^\alpha\bfxi^\beta \in \gr_w(\Dx).$$ We simply denote $\IN_w(0)=0$.

\begin{definition} Let $I$ be an ideal in $\Dx$. The \bem{characteristic ideal}
$$\cI(I) = \IN_{w}(I)\subset \gr_w(\Dx) = \bC[\bfx,\bfxi],$$
with $w = (0,e)$, i.e.: $w(x_i)=0$ and $w(\partial_i)=1$ for $i=1,\ldots,n$
defines the \bem{characteristic variety}
  $$\cV(\Dx/I) = \bV(\cI(I)) \subset \bC^{2n}.$$
\end{definition}

Here $\IN_w(I)$ stands for the ideal of $\gr_w(\Dx)$ generated by $\{\IN_w(Q), \vert\, Q\in I\}$, the set of the leading forms of the elements of $I$ with respect to the weight $w$. The set $\bV(J)$, for an ideal $J$ in $\bC[\bfx,\bfxi]$, is the variety in $\bC^{2n}$ defined by the ideal $J$.

If $M$ is a nonzero finitely generated $\Dx$--module, one can define, following \cite{Bernstein70}, its characteristic
variety $\cV(M)\subset \bC^{2n}$ by using any {\em good filtration} on $M$. The set $\cV(M)$ is an affine algebraic  subset of $\bC^{2n}$. In particular one can consider its {\em Krull dimension} $\dim(\cV(M))$.

To each irreducible component $C$ of $\cV(M)$ one can associate its {\em multiplicity} $m_C$ which is a positive integer number.
For the definition of the Krull dimension of an affine algebraic set and the multiplicity of its irreducible components one can consult \cite{Bourbaki-AC-Chap-8-9}.

\begin{theorem}[Bernstein's inequality; \cite{Bernstein70}] Let $M$ be a nonzero, finitely generated $\Dx$--module. Then
${n}\leq\dim(\cV(M))\leq 2n$. In particular, if $I \varsubsetneq \Dx$ is a left ideal, then ${n}\leq\dim(\cV(\Dx/I))\leq 2n$.
\end{theorem}

There is a number of software implementations of $D$-module algorithms:
{\tt kan/sm1}~\cite{KANwww},
{\tt risa/asir}~\cite{risa-asir-www},
{\tt dmod.lib} library~\cite{dmod-singular} in Singular~\cite{Singular},
and the {\tt Dmodules} package~\cite{DmodulesM2} in Macaulay2~\cite{M2}.
The last package hosts our implementations of the algorithms in this article.

\subsection{Holonomic $D$-modules}
A finitely generated $D_X$-module $M$ is called \bem{holonomic} if either $M=(0)$ or the dimension of its characteristic variety equals~${n}$ (i.e. $\dim(M)=n$). An ideal $I$ in $\Dx$ is called holonomic if either $I=\Dx$ or $\Dx/I$ is holonomic. 

The characteristic cycle of a holonomic $\Dx$--module $M$ is by definition the sum $$ \CChar(M) := \sum_{C {\mbox{ \footnotesize{irreducible component of }}} \cV(M)} m_C C$$ viewed as a cycle in the cotangent bundle $T^*\bC^n$. Here $m_C$ is the multiplicity of $C$ in $\cV(M)$.

Holonomic modules are particularly nice from the point of view of computations due to the following fact.
\begin{theorem}[Stafford \cite{StaffordModuleStructure}]
Every holonomic $\Dx$-module is cyclic.
\end{theorem}
Every holonomic module $M$ can be thought of as $M = \Dx\cdot Q$, where $Q$ is a cyclic generator, hence,  $M \cong \Dx/\Ann_\Dx(Q).$

One example of a holonomic module is $R=\bC[\bfx]$; it is generated by $1\in R$. Another is its localization $R_f=\bC[\bfx,f^{-1}]$, as stated by a theorem by J. Bernstein.

\begin{theorem}[\cite{Bernstein70, Bernstein}] Let $f\in \bC[\bfx]$ be a nonzero polynomial. Then the $\Dx$-module $R_f$ is holonomic. Moreover, for some negative integer $a$ the element $f^a\in R_f$ generates $R_f$ as a $\Dx$-module.
\end{theorem}

Computing the annihilator $\Ann_\Dx(f^a)$ is the main topic of this article.

By \cite{Bernstein70, Bernstein} the largest possible exponent $a = a(f)$ above is equal to the smallest integer root of the \bem{Bernstein-Sato polynomial} $b_f(s) \in \bQ[s]$, defined as the monic nonzero polynomial $b(s)$ of the smallest possible degree such that
$$
        b(s) f^s = Q(s,\bfx,\bfp)\cdot f^{s+1}, \mbox{ where } Q\in \Dx[s]:=\bC[s]\otimes_\bC \Dx.
$$

The idea of many localization algorithms (first formulated by Oaku~\cite{Oaku:local-b} in 1996) is simple to state:
\begin{enumerate}
      \item Compute  $\Ann(f^s)$ as an ideal of $\Dx[s]$.
      \item Compute $b_f(s)$ to determine $a=a(f)$.
      \item Specialization: $\left.\Ann(f^s)\right|_{s=a}$.
\end{enumerate}

This approach, for example, is implemented in the function {\tt Dlocalize} of the \D4M2\ package.
The bottlenecks of the algorithm are items 1 and 2 that require an expensive elimination via Gr\"obner bases in the Weyl algebra $\Dx$.

Step~2 can be avoided by using the following estimate
\begin{theorem}[\cite{Kashiwara-78}, see also \cite{Varchenko-82}, \cite{Saito-94}]\label{thm:lower-bound}
  $R_f$ is generated by $f^{-n+1}$.
\end{theorem}

As a corollary, for planar curves $f\in \bC[x,y]$, localization $R_f= \Dx\cdot f^{-1} = \Dx /\Ann (f^{-1})$.

\begin{remark}\label{rem:rat-fun-WA-annihilator}
  Note that one can easily compute the annihilator of a rational function in the the algebra of linear differential operators with \bem{rational function} coefficients $\bC(\bfx)[\bfp] = \bC(\bfx)\otimes_{\bC[\bfx]}\Dx$: the annihilator is a maximal ideal of $\bC(\bfx)[\bfp]$, e.g.,
  $$
  \Ann_{\bC(\bfx)[\bfp]}(f^a) = \bC(\bfx)\otimes \Ann_{\Dx}(f^a) =\<f\p_i - a\frac{\p f}{\p x_i} \,|\, i=1,\ldots,n \>.
  $$
\end{remark}

\section{Iterative algorithm}
Our main example will be the family of the so-called Reiffen's curves; see \cite{Reiffen-67}.

$$
f = f_{p,q} = x^p+y^q+xy^{q-1} = 0,
$$
where $p\geq 4, q\geq p+1$.

\begin{example} Computing the annihilator (this can be done via {\tt AnnFs} followed by specialization $s=-1$),
    \begin{eqnarray*}
      \Ann(f_{4,5}^{-1}) &=& \<4 x^{2} \p_x+5 x y \p_x+3 x y
      \p_y+4 y^{2} \p_y+16 x+20 y,\\
      &&16 x y^{2} \p_x+4y^{3} \p_x+12 y^{3} \p_y-125 x y \p_x-4 x^{2}\p_y+\\
      &&\ \ 5 x y \p_y-100 y^{2} \p_y+64 y^{2}-500 y,\\
      &&16 y^{3} \p_x^{2}-16 y^{3} \p_x \p_y+125 x y\p_x^{2}-
      35 x y \p_x \p_y+100 y^{2} \p_x\p_y+\\
      &&\ \ 12 x^{2} \p_y^{2}- 2 x y \p_y^{2}-24 y^{2}\p_y^{2}+112 x y \p_x-
      36 y^{2} \p_x+\\
      &&\ \ 84 y^{2}\p_y-930 x \p_x+625 y \p_x+26 x \p_y-\\
      &&\ \ 893 y\p_y+448 y-3720\>,
    \end{eqnarray*}
    we see that it is generated \bem{in order} $2$, i.e., generated by operators of order at most $2$.

    One interesting question is: What is the order of generation of this annihilator for a given $p$ and $q$?
\end{example}

\subsection{Iterative approach}

\begin{definition}\label{def:truncated-ann}
We call
  $$\Ann^{(d)}(f^a)= \left\< Q \in \Ann(f^a)\,|\, \ord Q \leq d \right\>$$
the $d$-th truncated annihilator of $f^a$.
\end{definition}

To compute $\Ann^{(d)}(f^{a})$ one may find the ${R}$-syzygy module $S_d$ for the vector of partial derivatives $(\bfp^\alpha\cdot f^{a})_\alpha$ where $|\alpha|\leq d$.

\begin{example} To find an element in $\Ann^{(1)}(f^{-1})$ for $f = x^2-y^3$ consider all partial derivatives of $f^{-1}$ of order at most 1:
$$\begin{array}{rccc}
  \p_x \cdot f^{-1} &=&-2x&f^{-2}\\
  \p_y \cdot f^{-1} &=&3y^2&f^{-2}\\
  1 \cdot f^{-1} &=&(x^2-y^3)&f^{-2}
\end{array}$$
For instance, ${3x}(-2x)+{2y}(3y^3)+{6}(x^2-y^3) = 0$, hence,
$$3x\p_x+2y\p_y+6 \in \Ann^{(1)}(f^{-1}).$$
\end{example}

In \D4M2\ the function for computing the truncated annihilator is called {\tt kOrderAnnFa}:

\begin{M2}
\begin{verbatim}
i1 : loadPackage "Dmodules";

i2 : R = QQ[x,y];

i3 : f = x^2-y^3;

i4 : A1 = kOrderAnnFa(1,f,-1)

                                 2             3      2       2
o4 = ideal (3x*dx + 2y*dy + 6, 3y dx + 2x*dy, y dy - x dy + 3y )

\end{verbatim}
\end{M2}

Note that {\tt kOrderAnnFa(d,f,-1)} would return the same ideal for all $d\in\bN$ as the annihilator for this particular example is generated in order~$1$.

\begin{remark}
  Suppose $S_{d-1}$, or in other words $\Ann^{(d-1)}(f^{a})$, is known. Then the computation of $S_d$, or the $d$-th truncated annihilator, can be optimized by computing the syzygies modulo $S_{d-1}$.

  In practice, however, it appears that the time of the syzygy computation at order $d$ dominates that for order $d-1$ for all $d$. In particular, in Experiment~\ref{exp:Reiffen} for the Reiffen curve of any degree the last step takes more time than all previous steps combined.
\end{remark}

\subsection{Stopping criterion}

The sequence of truncated annihilators stabilizes at some point:
$$
\Ann^{(1)}(f^{a}) \subseteq \Ann^{(2)}(f^{a}) \subseteq \cdots \subseteq \Ann^{(d)}(f^{a}) = \Ann^{(d+1)}(f^{a}) = \cdots
$$
\begin{definition}
   The smallest $d$ such that $\Ann^{(d)}(f^{a}) = \Ann(f^{a})$ is denoted by $\kappa(f^a)$ and is called the \bem{annihilator order} of $f^a$.
\end{definition}

Note that $\Ann^{(d)}(f^{a}) = \Ann^{(d-1)} (f^{a})$ does not imply $\Ann^{(d+1)} (f^{a}) = \Ann^{(d)} (f^{a})$.
For example, for $f=x$ and $a=3$ we have $Ann^{(1)} = Ann^{(2)} = Ann^{(3)} = \<x\p_x-3\>$ and $Ann^{(4)} = \<x\p_x-3, \p_x^4\>$.

\subsection{Annihilator order of a planar curve}

In view of Theorem~\ref{thm:lower-bound}, for a planar curve defined by $f\in R:=\bC[x,y]$ the localization $R_f$ is generated by $f^{-1}$. We have constructed an algorithm that computes truncated annihilators stopping exactly at the annihilator order of $f^{-1}$ thus producing the whole annihilator $\Ann(f^{-1})$ for a plane curve with at most one singular point,  that we can assume to be the origin.

\begin{algorithm} \label{alg:KappaAndAnn}
$(d,A) = $\KappaAndAnnihilator$(f)$

\begin{algorithmic}[1]
\REQUIRE $f \in R=\bC[x,y]$, a plane curve with at most one singular point, that we can assume to be the origin.
\ENSURE $d=\kappa(f^{-1})$, $A=\Ann(f^{-1})$.
\STATE $d \gets 0$.
\REPEAT
\STATE $d \gets d+1$.
\STATE $A \gets \Ann^{(d)}(f^{-1})$.
\STATE $\fa \gets \mbox{primary component of }\cI(A)$
corresponding to the origin. \label{step:primary-decomposition}
\UNTIL{$\deg \fa = (\mbox{multiplicity of }f \mbox{ at the origin}) - 1$}
\end{algorithmic}
\end{algorithm}
\begin{proof}[Proof of correctness]
We see $X=\bC^2$ as a complex manifold and consider $C$,
the analytic plane curve defined by the polynomial equation $f=0$. We assume $f$ to be a reduced polynomial. Denote \begin{itemize}
\item by $\cO_X$ the sheaf of holomorphic functions on $X$,
\item by $\cD_X$ the sheaf of linear partial differential operators on $X$ with holomorphic coefficients, and
\item by $\cO_X[*C]$ the sheaf of meromorphic functions on $X$ with poles on $C$.
\end{itemize}
By \cite{Kashiwara-78} the sheaf $\cO_X[*C]$ is a coherent holonomic $\cD_X$--module. If $p\in X\setminus C$ then $\cO_X[*C]_p$ is just $\cO_{X,p}$ and if the point $p\in C$ is a smooth point then the $\cD_{X,p}$--module  $\cO_X[*C]_p$ is isomorphic to the quotient
$$\frac{\cD_{X,(0,0)}}{\cD_{X,(0,0)} (x\partial_x+1,\partial_y)} $$
(for some choice of local coordinates $(x,y)$). So, in the neighborhood of a smooth point in $C$, the characteristic variety of $\cO_X[*C]$ is just the union of the conormal space to $C$, say $T^*_CX$ and the zero section  $T^*_XX$. Moreover, the multiplicities of each of these components in the corresponding characteristic cycle is 1.

It is well known that (see e.g. \cite[Section 6]{bri-merle-mais-1994}) 
in a sufficiently small neighborhood $U$ of a point $p\in C$, the characteristic cycle of $\cO_U[*C]$ is $$\CChar(\cO_U[*C]) = T^*_U U + \overline{T^*_{U\cap C\setminus \{p\}}U} + (m-1) T_{p}^* U$$ where $T^*_UU$ is the zero section of the cotangent  bundle $T^*U$, $\overline{T^*_{U\cap C\setminus \{p\}}U}$ stands for the closure in $T^*U$ of the conormal bundle to the smooth part of $U\cap C$, $T^*_pU$ is the conormal bundle to the point $p$ and $m$ is the multiplicity of the plane curve $C$ at the point $p$.

Let us denote $\calI^{} :=\calA nn^{}_{\cD_U}(f^{-1})$ and for any integer number $d\geq 1$,  $\calI^{(d)} :=\calA nn^{(d)}_{\cD_U}(f^{-1})$ considered as a sheaf of ideals in $\cD_U$.

As $f$ is a global equation for the plane curve $C\subset X$ we have $\calI^{(d)}_q = \cD_{U,q} Ann_{\Dx}(f^{-1})$ for any $q\in U$.

We have the following exact sequence of holonomic $\cD_U$-modules:
$$ 0 \rightarrow  \frac{\calI}{\calI^{(d)}} \rightarrow \frac{\cD_U}{\calI^{(d)}} \rightarrow  \frac{\cD_U}{\calI^{}} \rightarrow  0 $$
and then the following equality of characteristic cycles: $$ \CChar\left(\frac{\cD_U}{\calI^{(d)}}\right) = \CChar\left(\frac{\cD_U}{\calI^{}}\right) +  \CChar\left(\frac{\calI}{\calI^{(d)}}\right). $$

Denote by $m^{(d)}$ the multiplicity of the irreducible component $T^*_p U$ in $\CChar\left(\frac{\cD_U}{\calI^{(d)}}\right)$. As the ideals $\calI^{(d)}$ and $\calI^{}$ coincide on $U\cap C \setminus \{p\}$ the holonomic $\cD_U$--module $\frac{\calI}{\calI^{(d)}}$ is concentrated on $p$. Then $\CChar\left(\frac{\cD_U}{\calI^{(d)}}\right) = \CChar\left(\frac{\cD_U}{\calI^{}}\right)$ if and only $m^{(d)} = m-1$ (where $m$ is the multiplicity of the plane curve $C$ at $p\in C$). We also have $\calI^{(d)} = \calI$ if and only if $m^{(d)} = m-1$.

Once we have a finite set of polynomial equations for the characteristic ideal $\cI(Ann^{(d)}(f^{-1})),$ these equations define the characteristic cycle $\CChar\left(\frac{\cD_U}{\calI^{(d)}}\right)$ in the neighborhood $U$ of any point $p\in X$.
The multiplicity $m^{(d)}$ equals the degree of the primary ideal in the primary decomposition of $\cI(Ann^{(d)}(f^{-1}))$ corresponding to $T^*_pU$.
\end{proof}

\begin{remark}\label{rem:generic intersection}
The multiplicity $m^{(d)}$ can be also obtained as the multiplicity of the (unique) point $\tilde p\in T^*_p U$ in the intersection $\cI(Ann^{(d)}(f^{-1})) \cap L$ where $L$ is the defining ideal of a generic 2-plane in the 4-dimensional ambient space.
\end{remark}

Remark~\ref{rem:generic intersection} allows us to avoid the expensive primary decomposition step in the Algorithm~\ref{alg:KappaAndAnn}. Assume $p = 0$ is the only singularity; hence, let $U=X$ and let $S = \bC[x,y,\xi,\eta]$ be the ring of the characteristic ideal $\cI(Ann^{(d)}(f^{-1}))$. In our implementation we set $L = \<\xi-a,\eta-b\>$, $(a,b)\in\bC^2\setminus\{(0,0)\}$; in practice, $a$ and $b$ are taken to be small integers. The \bem{genericity} of the plane $L$ means that the point $\tilde p = (0,0,a,b)$ does not belong to $\overline{T^*_{U\cap C\setminus\{p\}}U}$. The following algorithm checks the genericity of the choice.

\begin{algorithm} \label{alg:CheckGenericity}
$g = $\CheckGenericity$(f,a,b)$

\begin{algorithmic}[1]
\REQUIRE $f \in R=\bC[x,y]$, a plane curve with the unique singular point at the origin.
\ENSURE $g$, a boolean value, whether the genericity condition described in Remark~\ref{rem:generic intersection} is satisfied.
\STATE Compute the syzygies $M \gets \{(u,v)\in R^2\,|\,uf_x+vf_y=0\}$.
\STATE $I = \< u\xi + v\eta\,|\,(u,v)\in S\> + \<f\>$.
\STATE Compute the saturation $J \gets I:\<x,y\>^\infty$.
\STATE $g \gets (0,0,a,b)\in \bV(J)$
\end{algorithmic}
\end{algorithm}

\begin{remark}
The choice of the plane $L = \<\xi-a,\eta-b\>$, $(a,b)\in\bC^2\setminus\{(0,0)\}$, is artificial but sufficient, since $\dim(T^*_p U \cap \bar{T^*_{U\cap C\setminus\{p\}}U})\leq 1$.
In fact, in Experiment~\ref{exp:Reiffen} the choice $(a,b)=(0,1)$ always works.
\end{remark}

\begin{remark}
In view of Remark~\ref{rem:generic intersection}, $m^{(d)}$ can be computed \bem{numerically}: one can pick a random $L$ and compute an approximation to the point of intersection $\tilde p$. Using this approximation, one may recover the multiplicity $m^{(d)}$ heuristically. To the best of our knowledge there is no known approach to certification of the correctness of such a computation.

Moreover, in our experiments the bottleneck is {\em not} the saturation step in Algorithm~\ref{alg:CheckGenericity}, but the computation of syzygies (described after Definition~\ref{def:truncated-ann}).
\end{remark}

\begin{remark}
  The algorithm can be generalized for arbitrary complex planar curves $C$ defined by a reduced polynomial $f$. As the number of singular points of $C$  is finite it is enough to apply our iterative algorithm at each singular  point of $C$ and compute $\kappa(f^{-1}_p)$ for each $p \in Sing(C)$. We finally have $\kappa(f^{-1})$ as the maximum of the  $\kappa(f^{-1}_p)$ for $p \in Sing(C)$. See also Section \ref{isolated-singularities} for a generalization of this algorithm to a more general situation.

\end{remark}

We have implemented Algorithm~\ref{alg:KappaAndAnn} modified according to the Remark~\ref{rem:generic intersection} and including the genericity check provided by Algorithm~\ref{alg:CheckGenericity} in the \D4M2\ package for Macaulay2.
The following example demonstrates that the new function {\tt kappaAnnF1PlanarCurve} produces the same result as the old method involving {\tt AnnFs}.
\begin{M2}
\begin{verbatim}
i2 : f = reiffen(7,8)

        7    8    7
o2 = x x  + x  + x
      1 2    2    1

o2 : QQ[x , x ]
         1   2

i3 : As = AnnFs f;

o3 : Ideal of QQ[x , x , dx , dx , s]
                  1   2    1    2

i4 : A = sub(As, {last gens ring As => -1});

i5 : (kappa,A') = kappaAnnF1PlanarCurve f;

i6 : kappa

o6 : 4

i7 : A == sub(A', ring A)

o7 = true
\end{verbatim}
\end{M2}

\begin{experiment}\label{exp:Reiffen}
We have computed the annihilator order for the Reiffen curve: $\kappa(f_{p,q}^{-1})$ seems to depend only on $p$. Starting with $p=4$ (the last computed value is for $p=21$), the sequence $\kappa(f_{p,*}^{-1})$ is
$$2,2,3,4,4,5,6,6,7,8,8,9,10,10,11,12,12,13,\ldots$$
We can compute $\kappa(f_{p,q}^{-1})$ for $q\leq 22$ in less than one day using our implementation of Algorithm~1. In contrast, an implementation of the general algorithm (the computation of $\Ann(f^s)$ followed by the specialization $s=-1$) took about a week for $p=17$ and has not produced a result in one month for $p=18$.
\end{experiment}

Another interesting observation is that the sequence of multiplicities $(m^{(d)})$ for $d=1,\ldots,\kappa(f_{p,q}^{-1})$ does not depend on $q$. Below we report the computed values $m^{(d)}$:
{
$$
\begin{array}{l|ccccccccccccc}
 d = & 1&2&3&4&5&6&7&8&9&10&11&12&13\\
\hline
p=4    & 4 & 3\\
p=5    & 6 & 4\\
p=6    & 8 & 6 & 5\\
p=7    & 10 & 8 & 7 & 6\\
p=8    & 12 & 10 & 9 & 7\\
p=9    & 14 & 12 & 11 & 9 & 8\\
p=10   & 16 & 14 & 13 & 11 & 10 & 9\\
p=11   & 18 & 16 & 15 & 13 & 12 & 10\\
p=12   & 20 & 18 & 17 & 15 & 14 & 12 & 11\\
p=13   & 22 & 20 & 19 & 17 & 16 & 14 & 13 & 12\\
p=14   & 24 & 22 & 21 & 19 & 18 & 16 & 15 & 13\\
p=15   & 26 & 24 & 23 & 21 & 20 & 18 & 17 & 15 & 14\\
p=16   & 28 & 26 & 25 & 23 & 22 & 20 & 19 & 17 & 16 & 15\\
p=17   & 30 & 28 & 27 & 25 & 24 & 22 & 21 & 19 & 18 & 16\\
p=18   & 32 & 30 & 29 & 27 & 26 & 24 & 23 & 21 & 20 & 18 & 17\\
p=19   & 34 & 32 & 31 & 29 & 28 & 26 & 25 & 23 & 22 & 20 & 19 & 18\\
p=20   & 36 & 34 & 33 & 31 & 30 & 28 & 27 & 25 & 24 & 22 & 21 & 19\\
p=21   & 38 & 36 & 35 & 33 & 32 & 30 & 29 & 27 & 26 & 24 & 23 & 21 & 20\\
      \end{array}
$$
}

\section{Discussion and open problems}

We have proposed a new algorithm computing a presentation, as a left $A_2$--module, of the localization $\bC[x,y]_f$ for nonzero polynomials $f=f(x,y)$ in two variables. The proof of the correctness of the algorithm is specific to dimension 2 and it seems difficult to generalize it to higher dimensions. Nevertheless, in the next subsection we describe an analogous algorithm for polynomials $f=f(x_1,\ldots,x_n)$, for any $n>0$, under the hypothesis that the hypersurface defined by $f$ in $\bC^n$ has only isolated singularities.

\subsection{Isolated hypersurface singularities}\label{isolated-singularities} The iterative algorithm described before can be also applied in a
more general situation. Let $n>0$ be an integer. Denote $X=\bC^n$,
by $\cO_X$ the sheaf of rigs of holomorphic functions on $X$ and by
$\cD_X$ the sheaf of rings of linear differential operators with
holomorphic coefficients. Assume that $f\in R=\bC[x_1,\ldots, x_n]$
is a non zero reduced polynomial defining a hypersurface $C \in X$
with an isolated singularity $p\in C$. Then the characteristic cycle
of $\cO_X[*C]$ has the following description (see e.g. \cite[Section 6]{bri-merle-mais-1994})
$$\CChar(\cO_X[*C]) = T^*_X X + \overline{T^*_{C\setminus \{p\}} X} + \nu T^*_p X$$ where $\nu$
is the Milnor number of a generic hyperplane section of $C$ through the point $p$. This formula coincides with the one we have considered before for $n=2$. The number $\nu$ can be computed in a effective way.
We can repeat the same ideas as in the proof of the correctness of Algorithm \ref{alg:KappaAndAnn}.  By comparison of the multiplicity $m^{(d)}$ with $\nu$, where $m^{(d)}$ is the multiplicity of the irreducible component $T^*_p X$ in $$\CChar\left(\frac{\cD_X}{\cD_X Ann_{A_n}(f^{-n+1})}\right),$$  we can compute the annihilator order of $f^{-n+1}$.

\subsection{Weyl closure}
\begin{definition}
The \bem{Weyl closure} of $I\subset \Dx$ is defined as $$\WeylCl(I) = K(\bfx)\otimes_{\bC[\bfx]} I \cap \Dx.$$
\end{definition}

Note that the problem of computing $\Ann(f^a)$ is a particular instance of the problem of computing the Weyl closure: according to \ref{rem:rat-fun-WA-annihilator} the annihilator of $f^{a}$ in $K(\bfx)\otimes_{\bC[\bfx]} \Dx$ is equal to $\<f\p_i + \frac{\p f}{\p x_i} \,|\, i=1,\ldots,n \>$. Therefore,
$$\Ann(f^{-1}) = \WeylCl(\<f\p_i + \frac{\p f}{\p x_i} \,|\, i=1,\ldots,n \>).$$

There exists an algorithm due to Tsai~\cite{HarryThesis} that computes the Weyl closure. However, the computation is highly nontrivial and boils down to computing localizations of $\Dx$-modules anyway. On the other hand, one may define the truncated Weyl closure of order $d$ analogously to the $d$-th truncated annihilator. Then, the sequence of truncated closures stabilizes at some degree; finding this degree or discovering a termination criterion for an iterative algorithm similar to Algorithm~\ref{alg:KappaAndAnn} would result in a much simpler procedure.

\section{Acknowledgements}
We thank the organizers of the Second CREST-SBM International Conference in Osaka, Japan, in 2010 where this project has been conceived. We thank L. Narv\'{a}ez and M. Barakat for their useful comments and suggestions. The second author thanks Institut Mittag-Leffler for hosting him in the Spring semester of 2011.

\ifx\form\JOURNAL
This work is partially supported by MTM2010-19336 and FEDER, FQM-333 and FQM-5849 Junta de Andaluc\'{\i}a
as well as by the US NSF grant DMS-0914802.
\fi

\end{document}